\documentclass[11pt,a4paper]{amsart}
\usepackage{amsmath, latexsym, amsfonts, amssymb,amsthm, amscd,epsfig,enumerate}
\usepackage{jipam}

\advance\hoffset -.75cm

\usepackage{amsfonts}
\usepackage{latexsym}
\usepackage{amsmath}
\usepackage{amssymb}

\newcommand{\E}{\mathbb E}

\newcommand{\R}{\mathbb R}

\newcommand{\N}{\mathbb N}
\newcommand{\pr}{\mathbb P}

\newcommand{\diff}{{\rm d}}

\newcommand{\I}{\mathbf 1}

\newcommand{\cF}{\mathcal{F}}

\newcommand{\calN}{\mathcal{N}}

\begin{document}

\title[An Inequality for correlated measurable functions]{An inequality for correlated measurable functions}


\author[F.~Zucca]{Fabio Zucca}
\address{F.~Zucca, Dipartimento di Matematica,
Politecnico di Milano,
Piazza Leonardo da Vinci 32, 20133 Milano, Italy.}
\email{fabio.zucca\@@polimi.it}
\urladdr{http://www1.mate.polimi.it/~zucca}


\keywords{integral inequalities, measure, cartesian product, ordered set}
\subjclass[2000]{26D15, 28A25}

\begin{abstract}
A classical inequality, which is known for families of monotone functions, is generalized to
a larger class of families of measurable functions. Moreover we characterize all the families
of functions for which the equality holds. We give two applications of this result, one
of them to a problem arising from
probability theory.
\end{abstract}

\maketitle



\section{Introduction}
\label{sec:intro}

The aim of this paper is to generalize an inequality,
originally due to Chebyshev and then rediscovered by Stein
in \cite{cf:Stein}. Usually this result is
stated for monotonic real functions: the classical inequality is
\[
(b-a)\int_a^bf(x)g(x)\diff x \ge \int_a^b f(x) \diff x \int_a^b g(x) \diff x
\]
where $f$ and $g$ are monotonic (in the same sense) real functions
(see for instance \cite{cf:Stein} and \cite{cf:Brualdi} for a
more general version).
If $a=b-1$ then this inequality has a probabilistic interpretation, namely
$\E[fg] -\E[f]\E[g] \ge 0$ (where $\E$ denotes the expectation), that is,
the covariance of $f$ and $g$ is nonnegative.

Our approach allows us to prove the inequality for
functions defined on a general measurable space,
hence we go beyond the usual ordered set $\R$.
More precisely, we prove an analogous result
for general families of measurable functions
that we call correlated functions (see Definition~\ref{def:corr}
for details). In particular we characterize all the families of functions for which
the equality holds.

Here is the outline of the paper.
In Section~\ref{sec:notation} we introduce the terminology and
the main tools needed in the sequel.
In particular
Sections~\ref{subsec:indord} and~\ref{subsec:indalg}
are devoted to the construction of an order relation
and a $\sigma$-algebra on a particular quotient space.
In Section~\ref{sec:main} we state and prove
our main result (Theorem~\ref{th:main}) which involves
$k$ correlated functions; the special case $k=2$ requires weaker
assumptions (see also Remark~\ref{rem:main}).
We give two applications of this inequality in Section~\ref{sec:final}:
the first one involves a particular class of power series, while the second one
comes from probability theory.



\section{Preliminaries and basic constructions}
\label{sec:notation}

We start from a very general setting.
Let us consider a set $X$, a partially
ordered space $(Y, \ge_Y)$ and a family $\calN = \{f_i\}_{i \in \Gamma}$ (where $\Gamma$ is an arbitrary set) of
functions in $Y^X$.
We consider the equivalence relation on $X$
\[
x \sim y \Longleftrightarrow f_i(x)=f_i(y), \ \forall i \in \Gamma
\]
and we denote by $X{/_\sim}$ the quotient space, by $[x]$ the
equivalence class of $x \in X$ and by $\pi$ the natural projection
of $X$ onto $X/_\sim$. Roughly speaking, by means of this
procedure, we identify points in $X$ which are not separated by
the family $\calN$.

To the family $\calN$ corresponds a natural counterpart
$\calN_\sim =\{\phi_{f_i}\}_{i \in \Gamma}$ of functions in $Y^{X/_\sim}$,
where, by definition,
$\phi_f([x]):=f(x)$, for all $x\in X$ and for every $f \in Y^X$ satisfying
\begin{equation}
\label{eq:count}
\forall x,y \in X:x \sim y \Longrightarrow f(x)=f(y)
\end{equation}
(this holds in particular for all the functions in $\calN$). It is
clear that the family $\calN_\sim$ separates the points of $X/_\sim$.
\hfill\break Given any function $g$ defined on $X/_\sim$
we denote by $\pi_g$ the function $g \circ \pi$; observe that
$\phi_{\pi_g}=g$ for all $g \in Y^{X/_\sim}$ and $\pi_{\phi_f}=f$
for every $f$ satisfying equation~\eqref{eq:count}. Clearly $g
\mapsto \pi_g$ is a bijection from $Y^{X/_\sim}$ onto the subset of
function in $Y^X$ satisfying
 equation~\eqref{eq:count}.
\hfill\break
Note that given $f,f_1 \in Y^X$ which satisfy equation~\eqref{eq:count} (resp.~$g,g_1 \in Y^{X/_\sim}$)
then $f \ge_Y f_1$ (resp.~$g \ge_Y  g_1$) implies
$\phi_f \ge \phi_{f_1}$ (resp.~$\pi_g \ge\pi_{g_1}$).

\subsection{Induced order}
\label{subsec:indord}

In order to prove Theorem~\ref{th:main} we
cannot take advantage, as in the classical formulation, of
an order relation on the set $X$.
Under some
reasonable assumptions (se Definition~\ref{def:corr} below) we can
transfer the order relation from $Y$ to $X/_\sim$ where we already defined a
family $\calN_\sim$ related to the original $\calN$. This will be enough for
our purposes.

\begin{definition}\label{def:corr}
The functions in $\calN$ are \sl{correlated} if, for all $i \in \Gamma$ and $x,y \in X$,
\begin{equation}
\label{eq:qm}
f_i(x)>_Y f_i(y) \Longrightarrow f_j(x) \ge_Y f_j(y), \, \forall j \in \Gamma.
\end{equation}
\end{definition}
We note that the definition above can be equivalently stated as follows:
for all $i,j \in \Gamma$ and $x \in X$,
\[
f_i^{-1}((-\infty, f_i(x))) \subseteq f_j^{-1}((-\infty, f_j(x)]).
\]
Besides, if $Y=\R$ with its natural order, then the functions in $\calN$ are correlated
if and only if for all $i,j \in \Gamma$ and $x,y \in X$,
\begin{equation}
\label{eq:qm1}
(f_i(x)-f_i(y))(f_j(x)-f_j(y)) \ge 0.
\end{equation}
In particular if $X$ is a totally ordered set and all the functions in $\calN$ are nondecreasing (or nonincreasing)
then they are correlated.

A family of correlated functions induces a natural
order relation on the quotient space $X/_\sim$.

\begin{lemma}
If the functions in $\calN$ are correlated then
the relation on $X/_\sim$
\[
[x] \ge_\sim [y] \Longleftrightarrow f_i(x) \ge_Y f_i(y), \, \forall i \in \Gamma
\]
is a partial order. If $(Y, \ge_Y)$ is a totally ordered
space then the same holds for $(X/_\sim, \ge_\sim)$.
Moreover $\calN_\sim$ is a family of nondecreasing functions
(hence they are correlated).
\end{lemma}
\begin{proof}
It is straightforward to show that $\ge_\sim$ is a well-defined partial order
(clearly it does not depend on the choice of $x$ (and $y$) within an equivalence class).
We prove that, if $\ge_Y$ is a total order, the same holds for $\ge_\sim$. Indeed
if $[x]\not = [y]$ then there exists $i \in \Gamma$ such that $f_i(x)\not = f_i(y)$;
suppose that $f_i(x) > f_i(y)$
then, by equation~\eqref{eq:qm},
$[x] >_\sim [y]$.
It is trivial to prove that $\phi_{f_i}$ is nondecreasing for every $i \in \Gamma$,
whence they are correlated since the space $(X/_\sim,\ge_\sim)$ is totally ordered.
\end{proof}

A subset $I$ of an ordered set, say $Y$, is called an interval if
and only if for all $x,y \in I$ and $z \in Y$ then $x \ge_Y z \ge_Y y$ implies
$z \in I$.
Note that given an interval $I \subseteq Y$ then $\phi_{f_i}^{-1}(I)$ is an interval of
$X/_\sim$ for every $i \in \Gamma$.

Given $x,y \in X$ such that $[x] \ge_\sim [y]$ we define the interval
$[[y],[x]):=\{[z] \in X/_\sim:  [y] \le [z] < [x]\}$;
the intervals $[[y],[x]]$, $([y],[x]]$ and $([y],[x])$ are
defined analogously.
In particular for any $x \in X$, we denote by $[[x],+\infty)$ and $(-\infty,[x]]$
the intervals $\{[y] \in X/_\sim: [y] \ge_\sim [x]\}$ and
$\{[y] \in X/_\sim: [x] \ge_\sim [y]\}$ respectively.

\subsection{Induced $\sigma$-algebra and measure}
\label{subsec:indalg}

This construction can be carried on under general assumptions. Let us consider
a measurable space with a positive measure $(X, \Sigma_X, \mu)$ and
an equivalence relation $\sim$ on $X$ such that for all $x \in X$ and $A \in \Sigma_X$,
\begin{equation}
\label{eq:equiv}
x \in A  \Longrightarrow  [x] \subseteq A.
\end{equation}
There is a natural way to construct
a $\sigma$-algebra on $X/_\sim$, namely define
\[
\Sigma_\sim :=\{\pi(A): A \in \Sigma_X\}
\]
where $\pi(A):=\{[x]:x \in A\}$. This is the largest $\sigma$-algebra
on $X/_\sim$ such that the projection map $\pi$ is measurable.
Observe that $A \mapsto \pi(A)$ is
a bijection from $\Sigma_X$ onto $\Sigma_\sim$.
It is natural to define a measure $\overline \mu:=\mu_\pi$
by
\[
\overline \mu(\pi(A))=\mu(A), \, \forall A \in \Sigma_X.
\]
It is well known that a function $g:X/_\sim \to \R$
is measurable if and only if $\pi_g$ 
is measurable.
%
%
%
%
Moreover $g$ is integrable (with respect to $\overline \mu$)
if
and only if $\pi_g$ is integrable (with respect to $\mu$) and
\begin{equation}
\label{eq:eqint}
\int_X \pi_g \diff \mu =
\int_{X/_\sim} g \diff \overline \mu.
\end{equation}
%
%
We say that a function $g$
is {\sl integrable} if at least
one of the integrals of the two nonnegative functions $g^+:=\max(g,0)$ and $g^-:=-\min(g,0)$ is finite;
hence the integral of $g$ can be unambiguously defined
as the difference of the two integrals (where
$\pm\infty+z:=\pm\infty$ for all $z\in \R$ and $0 \cdot \pm \infty:=0$).
This notion is sligthly weaker than the usual one: to remark
the difference, when the integrals of $g^+$ and $g^-$ are both finite
the function $g$ is called {\sl summable}.

It is a simple exercise to check that the equivalence relation
defined in Section~\ref{subsec:indord} satisfies equation~\eqref{eq:equiv}
if $\Sigma_X=\sigma(f_i:i \in \Gamma)$ (that is, $\Sigma_X$ is
the minimal $\sigma$-algebra such that all the functions in $\calN$
are measurable);
this equivalence relation along with its induced $\sigma$-algebra and
measure will play a key role in the next section.

\begin{remark}
\label{rem:GSfunct}
It is easy to show that if $h, r: X \mapsto \R$ are two integrable functions
such that the sum $\int_X h \diff \mu+\int_X r \diff \mu$ is not ambiguous
(i.e.~it is not true that $\int_X h \diff \mu=\pm \infty$ and $\int_X r \diff \mu=\mp \infty$)
then $h+r$ is integrable and
\begin{equation}
\label{eq:equality3}
\int_X (h+r)\diff \mu=\int_X h \diff \mu+\int_X r \diff \mu
\end{equation}
%
%
%
%
%
%
(both sides possibly being equal to $\pm \infty$). 
This will be useful in the proof
of Lemma~\ref{lem:main1}.
\end{remark}

\section{Main result}
\label{sec:main}

Throughout this section we consider a
measurable space with finite positive measure $(X,\Sigma_X,\mu)$
and a family of correlated functions $\calN=\{f_i\}_{i \in \Gamma}$,
where $\Sigma_X=\sigma(f_i:i\in \Gamma)$.
Let us consider $Y=\R$ with its natural order $\ge$. The equivalence relation $\sim$,
the (total) order $\ge_\sim$
and the space $(X/_\sim, \Sigma_\sim, \overline \mu)$ are
introduced according to Sections~\ref{subsec:indord} and~\ref{subsec:indalg}.
It is clear that  $\Sigma_\sim$ contains
the $\sigma$-algebra generated by the set of intervals
$\{\phi_{f_i}^{-1}(I): i \in \Gamma, I \subseteq \R \text{ is an interval}\}$.
More precisely it is easy to see that, by construction, all the intervals
of the totally ordered set $(X/_\sim,\ge_\sim)$ are measurable since
$\calN_\sim$ separates points.
%
%
%
%

The main result is the following.

\begin{theorem}
\label{th:main}
Let $\mu(X)<+\infty$.
\begin{enumerate}
\item
If $f$, $g$ are two integrable, $\mu$-a.e.~correlated functions such that
$fg$ is integrable then
\begin{equation}
\label{eq:main6}
\int_X fg \diff \mu \ge \int_X f \diff \mu \int_X g \diff \mu.
\end{equation}
Moreover, if $f$, $g$ are summable, then in the previous
equation the equality holds if and only if at least one of the functions
is $\mu$-a.e~constant.
\item
If $\{f_i\}_{i=1}^k$ be a family of measurable functions on $X$ which are
nonnegative and $\mu$-a.e.~correlated  then
\begin{equation}
\label{eq:main}
\mu(X)^{k-1} \int_X \prod_{i=1}^k f_i \diff \mu \ge
\prod_{i=1}^k \int_X f_i \diff \mu.
\end{equation}
Moreover if
$\int_X f_i \diff \mu \in (0,+\infty)$ for all $i=1,\ldots,k$,
then in the previous equation the equality holds if and only if at least
$k-1$ functions are $\mu$-a.e.~constant.
\end{enumerate}
\end{theorem}

Before proving this theorem, let us warm up with the following lemma;
though it will not be used in the proof of Theorem~\ref{th:main}, nevertheless
it sheds some light on the next step.

\begin{lemma}
\label{lem:main}
Let $\calN:=\{\{x_i(j)\}_{i \in \N}\}_{j=1}^k$ be a family of nonnegative and nondecreasing sequences
and $\{\mu_i\}_{i \in \N}$ be a family of strictly positive real numbers. If
$\sum_i \mu_i < +\infty$ then
\begin{equation}
\label{eq:main1}
\Big (\sum_i \mu_i \Big)^{k-1} \sum_i \prod_{j=1}^k x_i(j) \mu_i \ge \prod_{j=1}^k \sum_i x_i(j) \mu_i.
\end{equation}
Moreover if for every $j$ we have
$0<\sum_i x_i(j)<+\infty$ then
the equality holds if and only if at least $k-1$ sequences are constant.
\end{lemma}

\begin{proof}
We prove the first part of the claim for two finite sequences
$\{x_i\}_{i=1}^n$ and $\{y_i\}_{i=1}^n$, since the general case follows easily
by induction on $k$ and using the Monotone Convergence Theorem as $n$ tends to infinity.

It is easy to prove that
\begin{equation}
\label{eq:main2}
\sum_{i=1}^n \mu_i
\sum_{i=1}^n x_iy_i\mu_i
-\sum_{i=1}^n x_i\mu_i
\sum_{i=1}^n y_i\mu_i=\sum_{i,j:i \ge j} (x_i-x_j)(y_i-y_j) \mu_i \mu_j=\sum_{i,j:i > j} (x_i-x_j)(y_i-y_j) \mu_i \mu_j.
\end{equation}
Indeed
\[
\sum_{i=1}^n \mu_i
\sum_{i=1}^n x_iy_i\mu_i= \sum_{i,j:i > j} (x_iy_i+x_jy_j) \mu_i\mu_j +\sum_{i=1}^n x_iy_i \mu_i^2
\]
and
\[
\sum_{i=1}^n x_i\mu_i
\sum_{i=1}^n y_i\mu_i= \sum_{i,j:i >j} (x_iy_j+x_jy_i) \mu_i\mu_j +\sum_{i=1}^n x_iy_i \mu_i^2.
\]
This implies easily that
\[
\sum_{i=1}^n \mu_i
\sum_{i=1}^n x_iy_i\mu_i
-\sum_{i=1}^n x_i\mu_i
\sum_{i=1}^n y_i\mu_i \ge 0.
\]
%
%

If either at least $k-1$ sequences are constant or one sequence is equal to $0$, then we have an equality.
The same is true if $\sum_i x_i(j) \mu_i=+\infty$ for some $j$ and $\sum_i x_i(j) \mu_i>0$ for all $j$,
since both sides of equation~\eqref{eq:main1} are equal to $+\infty$.
On the other hand by using the first part
of the theorem and by taking the limit in equation~\eqref{eq:main2} as $n$ tends to infinity, for all $1 \le j_1<j_2\le k$,
\begin{equation}
\label{eq:main4}
\begin{split}
\Big (\sum_i \mu_i \Big)^{k-1} & \sum_i \prod_{j=1}^k x_i(j) \mu_i - \prod_{j=1}^k \sum_i x_i(j) \mu_i \\
& \ge
\Big (\sum_i \mu_i \Big) \sum_i x_i(j_1) x_i(j_2) \mu_i
\prod_{j \not = j_1,j_2} \sum_i x_i(j) \mu_i - \prod_{j=1}^k \sum_i x_i(j) \mu_i\\
&=\Big ( \prod_{j \not= j_1,j_2} \sum_i x_i(j) \mu_i \Big )
\sum_{i,i_1:i > i_1} (x_i(j_1)-x_{i_1}(j_1))(x_i(j_2)-x_{i_1}(j_2)) \mu_i \mu_{i_1}.
\end{split}
\end{equation}
If both $\{x_i(j_1)\}_i$ and $\{x_i(j_2)\}_i$ are nonconstant then there
exist $r<l$ and $r_1<l_1$ such that $x_r(j_1)<x_l(j_1)$ and
$x_{r_1}(j_2)<x_{l_1}(j_2)$. This implies
$x_{\max(l,l_1)}(j_1)-x_{\min(r,r_1)}(j_1)>0$ and
$x_{\max(l,l_1)}(j_2)-x_{\min(r,r_1)}(j_2)>0$,
thus the right hand side of equation~\eqref{eq:main4} is strictly positive
(just consider the summation over $\{i, i_1: i \ge \max(l,l_1), i_1 \le
\min(r,r_1)\}$) and
we have a strict inequality in equation~\eqref{eq:main1}.
\end{proof}

The proof of the previous lemma clearly
suggests a second lemma which will be needed
in the proof of Theorem~\ref{th:main}.
%
%

\begin{lemma}
\label{lem:main1}
Let $\calN:=\{f,g\}$ where $f,g:X \rightarrow \R$ are two summable functions such that
$fg$ is integrable (for instance if $f$ and $g$ are
$\mu$-a.e.~correlated).
If $ \mu(X)<+\infty$
then
\begin{equation}
\begin{split}
\label{eq:equality}
\mu(X)\int_{X} f(x)g(x)&\diff  \mu(x) =
\int_{X}f(x) \diff  \mu(x)
\int_{X}g(x) \diff  \mu(x)\\
&+
\frac12 \int_{X \times X} (f(x)-f(y))(g(x)-g(y))
\diff  \mu(x) \diff  \mu(y).
\end{split}
\end{equation}
\end{lemma}
%
%
%
%

\begin{proof}
Note that
\begin{equation}
\label{eq:equality1}
f(x)g(x)+f(y)g(y)=f(x)g(y)+f(y)g(x)
+
(f(x)-f(y))(g(x)-g(y));
\end{equation}
where $f(x)g(y)$ and $f(y)g(x)$ are summable  on $X \times X$, since $f,g$ are summable.
If we define $h(x,y):=f(x)g(y)+f(y)g(x)$ and $r(x,y):=(f(x)-f(y))(g(x)-g(y))$
then, according to Remark~\ref{rem:GSfunct}, we just need to prove that
$h$ and $r$ are integrable (since $h+r$ is integrable by hypothesis).

If $f$, $g$ are summable then, by equation~\eqref{eq:equality1},
$fg$ is integrable if and only if $(f(x)-f(y))(g(x)-g(y))$ is
integrable on $X \times X$ (since the sum of an summable
function and an integrable function is an integrable function)
and equation~\eqref{eq:equality} follows. Clearly if $f$ and $g$
are correlated then $(f(x)-f(y))(g(x)-g(y))$ is nonnegative thus
integrable.

\end{proof}

\begin{proof}[\sl{Proof of Theorem~\ref{th:main}}]
$ $\\
\begin{enumerate}
\item
By equation~\eqref{eq:eqint} it is enough to prove that
\[
\int_{X/_\sim} \phi_f \phi_g \diff \overline \mu \ge
\int_{X/_\sim} \phi_f \diff \overline \mu+
\int_{X/_\sim} \phi_g \diff \overline \mu.
\]
If $f$ and $g$ are summable then the claim follows from equation~\eqref{eq:equality}
of Lemma~\ref{lem:main1}. Otherwise, without loss of generality,
we may suppose that $\int_{X/_\sim} \phi_f \diff \overline \mu \equiv\int_{X} f \diff \mu=+\infty$.
%
%
If $\int_{X/_\sim} \phi_g \diff \overline \mu \equiv\int_{X} g \diff \mu < 0$
then there is nothing to prove. If $\int_X g \diff \mu  \ge 0$ then
either $g=0$ $\mu$-a.e.~, in this case both
sides of equation~\eqref{eq:main6} are equal to $0$, or
there exists $x \in X/_\sim$ such that $\overline \mu([x,+\infty))>0$ and
$\phi_f, \phi_g >0$ on $[x,+\infty)$ (since $\phi_f$ and $\phi_g$ are nondecreasing).
Clearly
$\int_{[x, +\infty)} \phi_f \diff \overline \mu =+\infty$
and
$\phi_f(y)\phi_g(y) \ge \phi_f(y)\phi_g(x)$ for all
$y \in [x, +\infty)$, hence both sides of equation~\eqref{eq:main6}
are equal to $+\infty$.

If one of the two functions is constant then the equality holds.
If $f$ and $g$ are nonconstant (that is, $\phi_f$ and $\phi_g$ are nonconstant)
then there exists $x_0,y_0 \in X/_\sim$ such that $x_0 >_\sim y_0$, $\phi_{f}(x_0)>\phi_{f}(y_0)$,
$\phi_{g}(x_0)>\phi_{g}(y_0)$,
$\overline \mu((-\infty, y_0])>0$ and $\overline \mu([x_0, +\infty))>0$
(this can be done as
in Lemma~\ref{lem:main1}).
Hence,
using equation~\eqref{eq:equality},
we have that,
\[
\begin{split}
\overline \mu(&X/_\sim)^{} \int_{X/_\sim} \phi_{f} \phi_{g} \diff \overline \mu -
\int_{X/_\sim} \phi_{f} \diff \overline \mu \int_{X/_\sim} \phi_{g} \diff \overline \mu\\
& \ge
\Big (
\int_{[x_0,+\infty)\times(-\infty,y_0]} (\phi_{f}(x)-\phi_{f}(y))(\phi_{g}(x)-\phi_{g}(y))
 \diff \overline \mu(x) \diff\overline \mu(y) \Big )\\
& \ge
\overline \mu((-\infty, y_0])\, \overline \mu([x_0, +\infty))
(\phi_{f}(x_0)-\phi_{f}(y_0))(\phi_{g}(x_0)-\phi_{g}(y_0))
>0.
\end{split}
\]
\item
Let us suppose that $f_i$ is  summable for all $i=1,\ldots,k$.
It is enough to prove that
\[
\overline \mu(X/_\sim)^{k-1} \int_{X/_\sim} \prod_{i=1}^k \phi_{f_i} \diff \overline \mu \ge
\prod_{i=1}^k \int_{X/_\sim} \phi_{f_i} \diff \overline \mu.
\]
In the previous part of the theorem,  we proved the claim
for two functions $\phi_f$ and $\phi_g$; as in Lemma~\ref{lem:main}, the general case follows by induction on $k$.


If at least two functions are nonconstant,
say $\phi_{f_1}$, $\phi_{f_2}$, then as before we may find
 $x_0,y_0 \in X/_\sim$ such that $x_0 >_\sim y_0$, $\phi_{f_1}(x_0)>\phi_{f_1}(y_0)$,
$\phi_{f_2}(x_0)>\phi_{f_2}(y_0)$,
$\overline \mu((-\infty, y_0])>0$ and $\overline \mu([x_0, +\infty))>0$
(this can be done as
in Lemma~\ref{lem:main1}).
By
applying the first part of the claim to the family (of $k-1$ functions)
$\phi_{f_1}\phi_{f_2}, \phi_{f_3}, \ldots, \phi_{f_k}$ (which are clearly still
correlated since they are nondecreasing) and
using equation~\eqref{eq:equality}
we have that,
\[
\begin{split}
\overline \mu(&X/_\sim)^{k-1} \int_{X/_\sim} \prod_{i=1}^k \phi_{f_i} \diff \overline \mu -
\prod_{i=1}^k \int_{X/_\sim} \phi_{f_i} \diff \overline \mu \\
&=
 \Big (\mu(X/_\sim) \int_{X/_\sim}  \phi_{f_1} \phi_{f_2} \diff \overline \mu -
\int_{X/_\sim} \phi_{f_1} \diff \overline \mu \cdot \int_{X/_\sim} \phi_{f_2} \diff \overline \mu \Big )
\prod_{i=3}^k \int_{X/_\sim} \phi_{f_i} \diff \overline \mu\\
& \ge
\Big (
\int_{[x_0,+\infty)\times(-\infty,y_0]} (\phi_{f_1}(x)-\phi_{f_1}(y))(\phi_{f_2}(x)-\phi_{f_2}(y))
 \diff \overline \mu(x) \diff\overline \mu(y) \Big )
\prod_{i=3}^k \int_{X/_\sim} \phi_{f_i} \diff \overline \mu\\
& \ge
\overline \mu((-\infty, y_0])\, \overline \mu([x_0, +\infty))
(\phi_{f_1}(x_0)-\phi_{f_1}(y_0))(\phi_{f_2}(x_0)-\phi_{f_2}(y_0))
\prod_{i=3}^k \int_{X/_\sim} \phi_{f_i} \diff \overline \mu
>0
\end{split}
\]
since $0<\int_{X/_\sim} \phi_{f_i} \diff \overline \mu < +\infty$ for all $i=1, \ldots, k$,
thus the second part of the claim is proved.
\end{enumerate}
\end{proof}

Note that if $\int_X f_i \diff \mu=+\infty$ for some $i$ and $\int_X f_j \diff \mu>0$ for all $j$
(otherwise both sides of equation~\eqref{eq:main} are equal to $0$) then
both sides of equation~\eqref{eq:main} are equal to $+\infty$; indeed
apply the first
part of the theorem to the family of correlated bounded functions $\{\min(f_i,n)\}_{i=1}^k$
(where $n\in\N$) and take the limit of both sides of equation~\eqref{eq:main}
as $n$ tends to $+\infty$.

\begin{remark}
\label{rem:main}
According to Theorem~\ref{th:main},
there is a difference between the case $k=2$ and $k >2$;
indeed in the latter case the inequality cannot be proved
for integrable (or even summable) $\mu$-a.e.~correlated functions
which are not nonnegative. Something happens
in the inductive process, namely if $\{f_i\}_{i=1}^k$ are correlated this may not be true for
$\{f_1f_2, f_3, \ldots, f_k\}$ (if the functions are not positive). Here is a counterexample:
take $X=[-1,1]$ endowed with the Lebesgue measure, $f_1(x)=f_2(x):=x\I_{[-1,0]}(x)$ and $f_i(x):=x-f_1(x)$ for all $i\ge3$.

Strictly speaking, Theorem~\ref{th:main} could be proved without the constructions
of Sections~\ref{subsec:indord} and~\ref{subsec:indalg}; one has just to use carefully
equation~\eqref{eq:qm1} and Lemma~\ref{lem:main1}. Our approach simplifies the proof of
Theorem~\ref{th:main} and gives a better understanding of the role of the correlation hypothesis
(compared to the usual monotonicity).

We finally observe that if we consider two integrable {\sl anticorrelated}
functions (meaning that $(f(x)-f(y))(g(x)-g(y)) \le 0$
for all $x,y \in X$) such that $fg$ is integrable then, clearly, we have
$\int_X fg \diff \mu \le \int_X f \diff \mu \int_X g \diff \mu$.
\end{remark}

\section{Final remarks and examples}\label{sec:final}

Let us apply Theorem~\ref{sec:main} to a class of power series.
We consider $f(z):=\sum_{n =1}^{+\infty} a_n z^n$ where
$\{a_n\}_n$ is a sequence of nonnegative real numbers and we suppose that
 $\{\rho^na_n\}$ is
nonincreasing (resp.~nondecreasing) for some $\rho$ such that $0 < \rho \le R$ (where
$R$ is the radius of convergence). Then the function
$z \mapsto (\rho-z)f(z)$ is a nonincreasing (resp.~nondecreasing)
on $[0,\rho)$.

Indeed if we suppose that $\{\rho^na_n\}$ is
nonincreasing then, for all $z, \gamma$ such that $0 \le z < \gamma <\rho$, we have
\[
\begin{split}
\sum_{n =1}^{+\infty} a_n z^n& = \sum_{n =1}^{+\infty} a_n \rho^n (z/\gamma)^n (\gamma/\rho)^n \\
&\ge \frac{\sum_{n =1}^{+\infty} a_n \gamma^n}{\sum_{n =1}^{+\infty} (\gamma/\rho)^n} \sum_{n =1}^{+\infty} (z/\rho)^n
=\sum_{n =1}^{+\infty} a_n \gamma^n \frac{\rho-\gamma}{\rho-z},
\end{split}
\]
where, in the first inequality, we applied Theorem~\ref{th:main} to the (correlated) functions
$f_1(n):=a_n\rho^n$ and $f_2(n):=(z/\gamma)^n$ defined on $\N$ endowed with the
measure $\mu(A):=\sum_{n \in A} (\gamma/\rho)^n$.
The case when $\{\rho^na_n\}$ is
nondecreasing is analogous (observe that now the functions $f_1$ and $f_2$ are anticorrelated).
If $z<\rho<R$ then $f_1$ and $f_2$ are nonconstant functions, hence the
function $z \mapsto (\rho-z)f(z)$
is strictly monotone.

\smallskip
We draw our second application application from probability theory.
To emphasize this, we denote the measure space by $(\Omega, \cF, \pr)$
and we speak of random variables and events
instead of measurable functions and measurable sets respectively.
We note that if $k=2$ then Theorem~\ref{th:main} says that
correlated variables have nonnegative covariance that is,
$\E[f_1f_2] - \E[f_1]\E[f_2]\ge 0$ (where $\E[f]:=\int_\Omega f \diff \pr$
is the usual expectation).

We call the (real) random variables $\{X_0,X_1,\ldots,X_k\}$
{\sl independent} if and only if, for every family of Borel sets
$\{A_0,A_1,\ldots,A_k\}$, we have $\pr(\cap_{i=0}^k \{X_i \in A_i\})=
\prod_{i=0}^k \pr(X_i \in A_i)$, where $\pr(X_i \in A_i)$
is shorthand for $\pr(\{\omega \in \Omega: X_i(\omega) \in A_i\})$.

In order to make a specific example, let us think of the variable $X_i$ ($i=1, \ldots, k$) as the
(random) time made by the $i$-th
contestant in an individual time trial bicycle race and let $X_0$ be our own (random) time; we suppose that
each contestant is unaware
of the results of the others (this is the independence hypothesis). If we
know the probability of winning a one-to-one race against each of our competitors we may be interested, for instance,
in estimating the probability of winning the race.
Such estimates are possible as
a consequence of Theorem~\ref{th:main}; indeed we have that
\begin{equation*}
\begin{split}
\pr(\cap_{i=1}^k \{X_i \ge X_0\})&\ge
\prod_{i=1}^k \pr(X_i \ge X_0)\\
\pr(\cap_{i=1}^k \{X_i \le X_0\})&\ge
\prod_{i=1}^k \pr(X_i \le X_0).\\
\end{split}
\end{equation*}
Thus the events $\{\{X_i \ge X_0\}\}_{i=1}^k$ (resp.~$\{\{X_i \le X_0\}\}_{i=1}^k$) are positively correlated (roughly speaking
this means that knowing that $\{X_1 \ge X_0\}$ makes, for instance, the event $\{X_2 \ge X_0\}$ more likely than before).

The proof of these inequalities is straightforward. If we define $\mu(A):=\pr(X_0 \in A)$ for all Borel sets $A \subseteq \R$,
then, according
to Fubini's Theorem,
\[
\begin{split}
\pr(X_i \ge X_0)=\int_\R \pr(X_i \ge t) \diff \mu(t), & \quad  \pr(\cap_{i=1}^k \{X_i \ge X_0\})=\int_\R \prod_{i=1}^k \pr(X_i \ge t)
\diff \mu(t)\\
\pr(X_i \le X_0)=\int_\R \pr(X_i \le t) \diff \mu(t), & \quad \pr(\cap_{i=1}^k \{X_i \le X_0\})=\int_\R \prod_{i=1}^k \pr(X_i \le t)
\diff \mu(t).\\
\end{split}
\]
Indeed
\[
\pr(X_i \ge X_0)=\int_{\{(s,t)\in \R^2:s \ge t\}} \diff \nu(s) \diff \mu(t)=\int_\R \int_{[t,+\infty)} \diff \nu(s) \diff \mu(t)=
\int_\R \pr(X_i \ge t) \diff \mu(t)
\]
where $\nu(A):=\pr(X_i \in A)$ for all borel sets $A \subseteq \R$ and the first equality holds since $X_i$ and $X_0$ are independent.
The remaining cases are analogous.
Note that $\{\pr(X_i \ge t)\}_{i=1}^k$ and $\{\pr(X_i \le t)\}_{i=1}^k$ are both
families of monotone (thus correlated) functions; Theorem~\ref{th:main} yields the claim.
This example can be easily extended to a more interesting case: namely when
$\{X_1, \ldots, X_k\}$ have identical laws and are independent conditioned to $X_0$ (see Chapters 4 and 6 of \cite{cf:Bill} for details).
In this case one can prove that
\[
\pr(\cap_{i=1}^k \{X_i \in A\})\ge
\prod_{i=1}^k \pr(X_i \in A), \qquad \forall A \subseteq \R \text{ Borel set}.
\]
The proof makes use of Theorem~\ref{th:main} in its full generality but
this example exceeds the purpose of this paper.

\section*{Acknowledgments}
The author thanks S.~Mortola for useful discussions.


\end{document}